\documentclass[10pt]{article}  

\usepackage{amsmath,amsthm,amssymb}
\usepackage{graphicx}

\newcommand\valpha{\vec \alpha}

\newcommand\vbeta{\vec \beta}
\newcommand\veca{\vec a}

\newtheorem{theorem}{Theorem}[section]
\newtheorem{lemma}[theorem]{Lemma}

\theoremstyle{definition}
\newtheorem{definition}{Definition}[section]

\theoremstyle{remark}
\newtheorem{remark}{Remark}

\newcommand{\bi}{\begin{itemize}} 
\newcommand{\ei}{\end{itemize}}
\newcommand{\ben}{\begin{enumerate}}
\newcommand{\een}{\end{enumerate}}
\newcommand{\be}{\begin{equation}}
\newcommand{\ee}{\end{equation}}
\newcommand{\bea}{\begin{eqnarray}} 
\newcommand{\eea}{\end{eqnarray}}
\newcommand{\ba}{\begin{align}} 
\newcommand{\ea}{\end{align}}
\newcommand{\bse}{\begin{subequations}} 
\newcommand{\ese}{\end{subequations}}
\newcommand{\bc}{\begin{center}}
\newcommand{\ec}{\end{center}}
\newcommand{\bfi}{\begin{figure}}
\newcommand{\efi}{\end{figure}}

\newcommand{\tbox}[1]{{\mbox{\tiny #1}}}  
\newcommand{\mbf}[1]{{\mathbf #1}}
\newcommand{\half}{\mbox{\small $\frac{1}{2}$}}

\newcommand{\fref}[1]{Fig.~\ref{#1}}          

\newcommand{\pA}{{{\mathcal A}^+}}       
\newcommand{\bigO}{{\mathcal O}}
\newcommand{\xx}{\mbf{x}}
\newcommand{\yy}{\mbf{y}}
\newcommand{\ui}{u^\tbox{inc}}
\newcommand{\us}{u^\tbox{sc}}
\newcommand{\GG}{{\Gamma}}         
\newcommand{\al}{\alpha}     

\def\keywords{\vspace{.5em}
{\textit{Keywords}:\,\relax%
}}

\title{A fast and robust solver for the scattering from a layered periodic structure containing multi-particle inclusions}

\author{Jun Lai\thanks{Courant Institute of Mathematical Sciences, New York University, NY 10012. Email: \textit{lai@cims.nyu.edu}}, Motoki Kobayashi\thanks{Courant Institute of Mathematical Sciences, New York University, NY 10012. Email: \textit{motokia.kobayashi@gmail} \textit{.com} }, and Alex Barnett\thanks{Department
of Mathematics, Dartmouth College, Hanover, NH 03755. Email: \textit{ahb@math.dartmouth.edu}}}

\date{\today}
\begin{document}
\maketitle
\begin{abstract}  
We present a solver for plane wave scattering from a
periodic dielectric grating 
with a large number $M$ of 
inclusions lying in each period of its middle layer.
Such composite material
geometries have a growing role in modern photonic devices and solar cells.
%
The high-order scheme is based on boundary integral equations,
and achieves many digits of accuracy with ease.
The usual way to periodize the integral equation---%
via the quasi-periodic Green's function---fails at Wood's anomalies.
We instead use the free-space Green's kernel for the near field, add
auxiliary basis functions for the far field,
and enforce periodicity in an expanded
linear system; this is robust for all parameters.
Inverting the periodic and layer unknowns,
we are left with a square linear system involving only the inclusion scattering
coefficients.
Preconditioning by the single-inclusion scattering matrix,
this is solved iteratively in $\bigO(M)$ time
using a fast matrix-vector product.
Numerical experiments show that
a diffraction grating containing
$M=1000$ inclusions per period can be solved to 9-digit accuracy
in under 5 minutes on a laptop.
\end{abstract}


\section{Introduction}
\label{Intro}

The modeling and design of periodic dielectric structures plays a central role in modern optics.
Tools such as diffraction gratings, photonic crystals, meta-materials, plasmonics, and other micro-scale structures,
are becoming key to efficient devices, including
lasers, sensors, anti-reflective surfaces and absorbers \cite{DCD1},
and solar cells \cite{atwater}.
For instance, in thin-film solar cell design \cite{santbergen,Krishnan:14}
the use of periodic structures, and nanoparticle inclusions, in
ordered or disordered composites, enhances absorption.
One then seeks a grating structure with a
specific arrangement of inclusions that maximizes absorption.
Other optimization problems include the design of photonic crystal
lenses \cite{optimlens}.
Related is the inverse problem
of inferring a structure from measurements \cite{scatterometry,Bao:14}.
Such tasks demand a large number of solutions of the direct
(forward) scattering problem.
Similar periodic and multi-particle wave scattering problems
arise in acoustics and elastodynamics,
and in general whenever a {\em super-cell} is used to approximate
the response of a random composite material (e.g. \cite{meisels07}).
Such considerations have spurred the development of
efficient methods for solving Helmholtz
and Maxwell frequency-domain boundary value problems in periodic geometries
\cite{DCD1,Bao:95,brunohaslam09,Barnett20106898,BG2011,GHL2014,delourme14,mlqp}.
High accuracy can be challenging to achieve
due to guided modes, resonances, and extreme parameter sensitivity.



Therefore, in this paper we consider the monochromatic
scattering from a layered periodic structure containing a large number $M$ of inclusions (``particles'') at given locations, as in a (generalized)
photonic crystal.
As shown in \fref{figure1}, the structure is periodic in the $x$ direction, layered in the $y$ direction and invariant along the $z$ direction.
Because of the two-dimensional (2D) geometry, there exist two fundamental polarizations in the electromagnetic scattering: transverse magnetic (TM) where the magnetic field is transverse to the $(x,y)$ plane, and transverse electric (TE) where the electric field is transverse to the $(x,y)$ plane. We will focus on TM polarization, noting that our technique applies to TE polarization without any essential difficulty.

The grating scattering problem has been mathematically very well studied.
It has been proved that for an arbitrary periodic dielectric and incident angle
the problem has a unique solution for all frequencies with the possible exception of a countable set of resonances (singular frequencies \cite{BoStar1994})
at which the solution is not unique.
Such physical resonances are not to be confused with {\em Wood's anomalies}
(for the definition see the next section),
which are frequencies where at least one of the Bragg diffraction orders
points along the grating, i.e.\ in the $x$ direction.
A Wood's anomaly does not prevent the solution from being unique,
although it does cause arbitrarily large sensitivity with respect to
the incident wave angle or frequency \cite{Linton2007165},
and also causes problems with certain integral equation methods
\cite{BG2011}.
One of the advantages of our scheme is that it
is applicable and accurate at or near Wood's anomalies, without
any modifications.

There exists a wide range of numerical methods for periodic diffraction,
including boundary integral equations \cite{arens06,brunohaslam09,BG2011,GHL2014,delourme14,mlqp}, finite element methods \cite{Baofinite,BaoYang1999}, Fourier expansion based methods \cite{RCWAmoharam}, and continuation methods \cite{Bruno:93}. In the time domain, the finite difference scheme has been discussed in \cite{HolSte2002}. The advantages of the integral approach over finite elements and finite differences are that it reduces the dimension by one (vastly reducing the number of unknowns), and achieves high-order accuracy with appropriate surface quadratures. However, the resulting linear system is often dense, making a naive matrix-vector product expensive when the number of unknowns is large. In this paper, we will reduce this cost via the fast multipole method (FMM) \cite{GR1987}.
 
More specifically, we propose an integral approach based on the free space Green's function; this bypasses the considerable complexities of computing
the periodic Green's function \cite{kurkcu,delourme14}. We split the representation of the scattered field in the grating structure into near field and far field components. The near field is represented by standard free-space Helmholtz single- and double-layer potentials on the material interfaces, while the far field is taken care by a {\em local expansion} (Fourier--Bessel or $J$ expansion) whose coefficients are fixed by enforcing the periodic boundary
condition explicitly in the linear system.
This builds upon recent ideas of the last author and co-workers
\cite{Barnett20106898,BG2011,mlqp}.

Solving for discretized layer densities on each of the $M$ inclusion boundaries would
introduce an unnecessarily large number of unknowns.
Hence, following \cite{GG2013,Lai:14},
we precompute the inclusion {\em scattering matrices},
then treat the set of outgoing scattering coefficients
as a reduced set of unknowns.
When particles are sub-wavelength, and not extremely close to each other,
this is highly accurate
with only 20 or so unknowns per particle \cite{Lai:14}.
The full rectangular linear system then couples these to the grating interface densities and periodizing $J$-expansion coefficients.
By eliminating the last two (via a Schur complement and pseudoinverse)
we are left with a square linear system for the particle scattering
coefficients, which we precondition with a block-diagonal matrix and then
solve via GMRES with FMM acceleration, with effort scaling linearly in $M$.
The result is a robust, efficient, high-order accurate solver
that we expect to be useful
for design and optimization problems for periodic photonic devices.


The outline of the paper is as follows. Section \ref{Form} gives the mathematical formulation of the periodic problem. Section \ref{scat_peri} proposes the integral approach for the scattering from a periodic structure without particle inclusions, based on the free space Green's function. Section \ref{scat_multi} reviews classical multi-particle scattering and discusses the evaluation of the scattering matrix. The quasi-periodizing scheme combining all the above techniques is given in Section \ref{Quasi_scheme}, and numerical experiments are shown in Section \ref{Numerical}. We draw conclusions in Section \ref{Conc}.

\section{Problem formulation}
\label{Form}

\begin{figure}[tbp]  
\centering
\includegraphics[width=10cm,height=6cm]{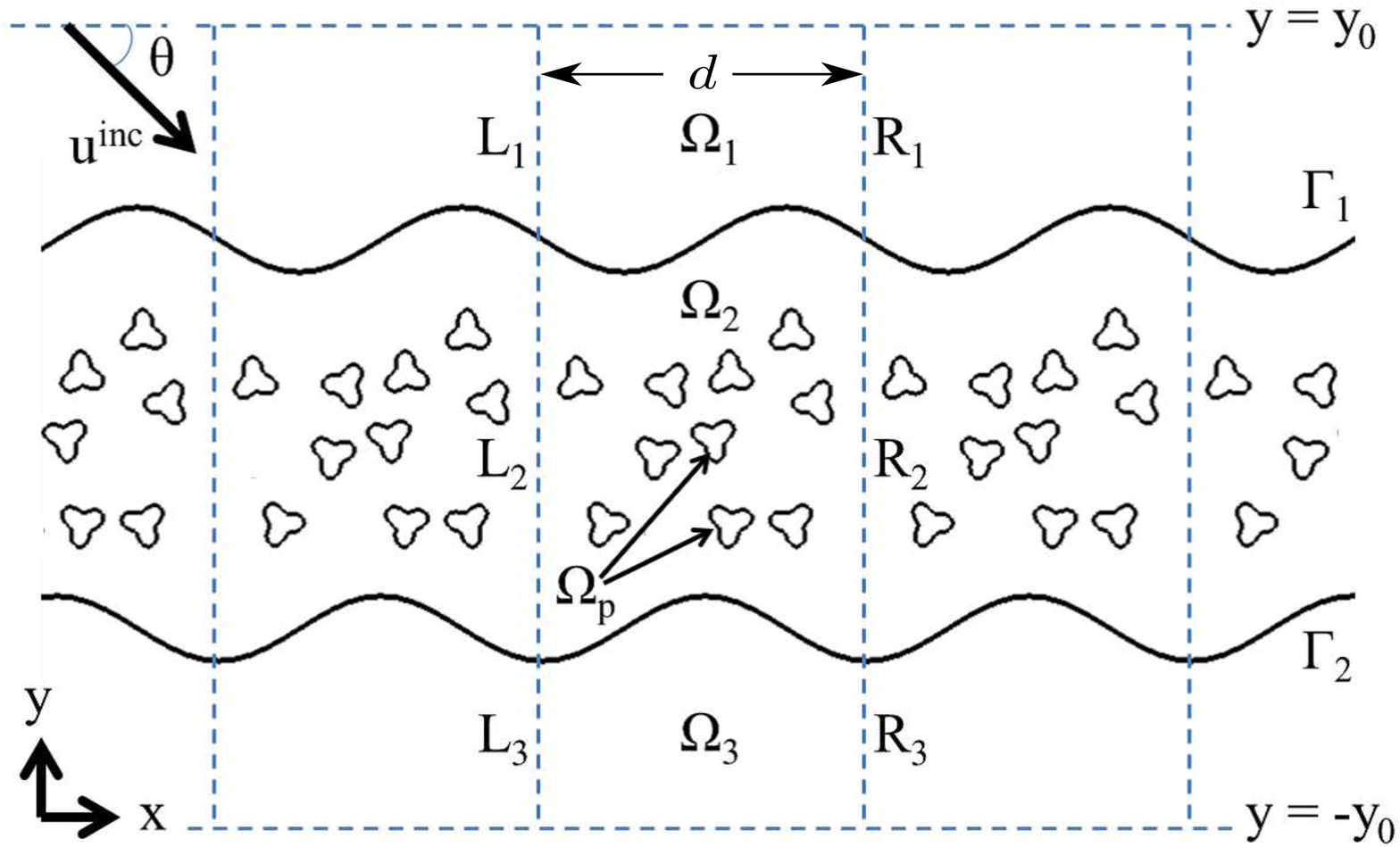}
\caption{A 2D grating scattering geometry: a plane wave incident on a three-layered periodic structure with periodicity $d$. A large number of identical
dielectric obstacles are embedded in the middle layer. We use $\Omega_p$ to denote the set of all these particles and $\Omega_p^j$ to denote the $j$th particle.
The vertical dotted lines indicate the unit cell walls $x= d/2\pm ld$, $l=1,2,3,\cdots$, while the
top and bottom dotted lines indicate the fictitious boundaries at $\pm y_0$. The three layers are denoted by $\Omega_1$, $\Omega_2$ and $\Omega_3$.}
\label{figure1}
\end{figure}

Consider the plane-wave incident time harmonic scattering (with time dependence $e^{-i\omega t}$) from a 2D periodic (or grating) structure with period $d$.
As shown in Figure \ref{figure1}, the unit cell $\Omega=[-d/2,d/2] \times \mathbb{R}$ consists of three layers, denoted by $\Omega_1$, $\Omega_2$ and $\Omega_3$. Let $\Gamma_1$ and $\Gamma_2$ denote the two smooth interfaces separating the layers.
The left and right boundaries of $\Omega_j$ are denoted by $L_j$ and $R_j$, $j = 1,2,3$. Assume the permittivity $\varepsilon$ is given as $\varepsilon_1$, $\varepsilon_2$ and $\varepsilon_3$ in the three layers respectively. A large number $M$ of particles, collectively denoted by $\Omega_p$, with the same permittivity $\varepsilon_p$, are located inside $\Omega_2$. The permeability $\mu$ is assumed to be constant everywhere.

For TM polarization, in which case the total electric field is $E(x,y) = (0,0,u)$, the full time harmonic Maxwell equations
\begin{align*}
\begin{cases}
\nabla\times E &= i\omega \mu H \\
\nabla\times H &= -i\omega\varepsilon E  
\end{cases}
\end{align*}
are reduced to the Helmholtz equation:
\begin{equation}
\label{eq1}
\Delta u +k(\xx)^2 u = 0~,
\end{equation}
where $\xx := (x,y)$, and
where the wavenumber $k$ takes one of four values,
\begin{align}
k(\xx) =\begin{cases}
k_1 := \omega\sqrt{\mu\varepsilon_1}, \qquad \xx \in \Omega_1\\
k_2 := \omega\sqrt{\mu\varepsilon_2}, \qquad \xx \in \Omega_2 \backslash \overline{\Omega_p}\\
k_3 := \omega\sqrt{\mu\varepsilon_3}, \qquad \xx\in\Omega_3\\
k_p := \omega\sqrt{\mu\varepsilon_p}, \qquad \xx\in \Omega_p
\end{cases}
\end{align}

In the usual setting of scattering theory, the full wave is $u = \ui+\us$, where $\ui$ is the incident wave and $\us$ is the resulting wave scattered from the periodic structure.
The incident wave is a plane wave $\ui(x,y) = e^{ik_1\cos\theta x+ik_1\sin\theta y}$ in $\Omega_1$, and $\ui=0$ elsewhere, with $\theta$ the incident angle. Since the wave is propagating in different layers, the continuity condition along various interfaces in TM polarization with constant permeability is:
\be
[u] = 0 ~, \hspace{1in}
\left[\frac{\partial u}{\partial n} \right] = 0 ~,
\label{cont}
\ee
where $[\cdot]$ denotes the jump of a function across the interface,
$\partial /{\partial n}$ is the normal derivative, and $u$ is the total field in each layer \cite{Cot}.

We use the term {\em quasi-periodic} if a function (such as $u$) satisfies
\be
u(x+d,y) = e^{i\kappa d} u(x,y) \qquad \mbox{for all } (x,y)~,
\label{quasi}
\ee
where $\kappa=k_1\cos\theta$ is the incident horizontal wavevector.
The factor $$\alpha:= e^{i\kappa d}$$ is the Bloch phase associated with translation by one period.
Since $\ui$ is quasi-periodic, we seek a scattered
wave with this same symmetry, hence \eqref{quasi} also holds for the full wave
\cite{Petit1980}.
Restricting to a single unit cell, with left and right walls $L$ and $R$ respectively, we have matching conditions
\bea
\al u_L - u_R &=& 0
\label{qp}
\\
\al \frac{\partial u}{\partial n}\biggr|_L - \frac{\partial u}{\partial n}\biggr|_R &=& 0
\label{qpn}
\eea
which are in fact equivalent to \eqref{quasi} \cite[Sec.~3]{BoStar1994}.

Finally, $\us$ must satisfy a radiation condition. Let $y_0$ be sufficiently large such that $\Omega_2$ lies between the lines $\Gamma_u: = \{y = y_0\}$ and $\Gamma_d := \{y = -y_0\}$ (see Figure \ref{figure1}). Define $\kappa_n=\kappa+2\pi n /d$, $n\in \mathbb{Z}$, and $k_{j,n} = +\sqrt{k_j^2-\kappa_n^2}$, where $j = 1,3$ and the sign of the square-root is taken as positive real or positive imaginary. Then the plane wave with wavevector $(\kappa_n,k_{j,n})$ is quasi-periodic for each integer $n$, and satisfies the Helmholtz equation at frequency $\omega$. The radiation condition on $\us$ is expressed by
uniformly convergent {\em outgoing or decaying} Rayleigh--Bloch expansions:
\begin{align}
\us(x,y) &= \sum_{n\in \mathbb{Z}}c_ne^{i\kappa_nx}e^{ik_{1,n}(y-y0)}, \qquad \mbox{ for } y>y_0, x\in \mathbb{R}, \label{rad1}\\
\us(x,y) &= \sum_{n\in \mathbb{Z}}d_ne^{i\kappa_nx}e^{ik_{3,n}(-y-y0)}, \qquad \mbox{ for } y<-y_0, x\in \mathbb{R}. \label{rad2}
\end{align}

The complex coefficients $c_n$, $d_n$, for $n$ such that $|\kappa_n|\le k_j$ (propagating waves), are the Bragg diffraction amplitudes at the grating orders. For all other $n$ these give evanescent components which do not contribute to the far field.
A Wood's anomaly occurs if, for some $n$, $\kappa_n=k_1$,
thus $k_{1,n}=0$
(upper half-space), or if $\kappa_n=k_3$
thus $k_{3,n}=0$
(lower half-space).
The radiation conditions ensure that $\us$ is outgoing except at a Wood's anomaly, when the $n$th Rayleigh--Bloch mode in \eqref{rad1} or \eqref{rad2} is
constant in $y$ (a horizontally traveling plane wave).
It is also possible to have a double Wood's anomaly, when $\kappa_n = \kappa_m = k_1$ for a pair of integers $n$, $m$.

The following theorem \cite{Baofinite,BaoDo2000,BoStar1994} describes the
well-posedness of the boundary value problem.
\begin{theorem}
Fixing an incident angle $\theta$, there
exists a unique solution $u$ to \eqref{eq1}--\eqref{cont}
and \eqref{qp}--\eqref{rad2} for all but a discrete set of $\omega$.
\end{theorem}
Henceforth we will assume parameter values $(\theta,\omega)$ for which
the solution is unique.

To summarize, we are interested in the solution of \eqref{eq1} together with the continuity condition \eqref{cont},
the unit-cell quasi-periodicity \eqref{qp}--\eqref{qpn},
and the radiation condition \eqref{rad1}--\eqref{rad2} satisfied by $\us$.
In the next section we discuss the solution for the layered periodic structure without inclusions, via the integral equation approach.

\begin{remark}
From now on we slightly abuse the notation $\Omega_1$, $\Omega_2$ and $\Omega_3$ introduced, by replacing them by their truncated versions. In other words, $\Omega_1$ has boundaries $L_1$, $\Gamma_1$, $R_1$, and $\Gamma_u$,
while $\Omega_3$ has boundaries $L_3$, $\Gamma_d$, $R_3$, and $\Gamma_2$.
The unit cell will refer to $\Omega =\Omega_1\cup \Omega_2\cup \Omega_3 = [-d/2,d/2] \times [-y_0,y_0]$.
The artificial upper and lower boundaries $\Gamma_u$ and $\Gamma_d$ are also called transparent boundaries in some literature \cite{Bao:95}. 
\end{remark}

\section{Robust solution for empty periodic layered structure}
\label{scat_peri}

The standard approach to convert a non-periodic integral equation
formulation into the correct periodic one is to use 
the quasi-periodic Green's function \cite{linton98}, defined 
at wavenumber $k_j$ as
\begin{equation}
\label{qGreen}
G^\tbox{qp}_{k_j}(\mathbf{x},\mathbf{y}) := \sum_{l \in \mathbb{Z}}\al^lG_{k_j}(\mathbf{x},\mathbf{y}+(ld,0))
~,
\end{equation}
where $\xx\in\Omega_j$ is the target and $\yy\in\Omega_j$ the source point,
and $G_{k}$ is the free space Green's function at wavenumber $k$, i.e.\ $G_{k}(\mathbf{x},\mathbf{y}) = 
\frac{i}{4}H_0^{(1)}(k\|\mathbf{x}-\mathbf{y}\|)$, 
where $H_0^{(1)}$ is the outgoing Hankel function of order zero,
and $\|\cdot\|$ is the Euclidean norm.
\eqref{qGreen} is well-defined away from
Wood's anomalies, and has been successfully applied in many grating problems \cite{arens06,brunohaslam09,GHL2014,Petit1980}. However, it has two major practical
drawbacks: i) it is expensive to evaluate (requiring either series
acceleration, or lattice sums \cite{linton98}),
and ii) it blows up (with an inverse-square-root singularity) at
Wood's anomalies, causing a purely numerical breakdown in the
solution of a what remains a well-posed problem.
%
The key idea is that \eqref{qGreen} can be rewritten as
\begin{equation}
\label{qGreen2}
G^\tbox{qp}_{k_j}(\mathbf{x},\mathbf{y}) = \sum_{l=-P}^{P}\al^l G_{k_j}(\mathbf{x},\mathbf{y}+(ld,0)) + \sum_{n\in\mathbb{Z}}a_n J_n(k_j\|\mathbf{x}-\mathbf{x_0}\|)
e^{i n \theta_{\xx-\xx_0}}
\end{equation}
where $P$ is a positive integer, $\xx_0\in\Omega_j$ is a fixed origin, $\theta_{\xx}$ is the angle of a vector $\xx$,
$J_n$ is the Bessel function of order $n$, and the {\em lattice sum}
coefficients $\{a_n\}$ can be found by Graf's addition theorem \cite[10.23(ii)]{Hand2010}.
The second term accounts for the smooth field due to the
infinite set of far sources $l<-P$ and $l>P$ from \eqref{qGreen}.
The $2P+1$ direct terms account for the near field.
\begin{remark}
The sum in \eqref{qGreen2} converges (and exponentially fast)
if and only if the target is closer to the
origin than the nearest ``far'' source, i.e.\
$\|\xx-\xx_0\| < \min_{l\in\mathbb{Z},l\notin [-P,P]} \|\yy + (ld,0) - \xx_0\|$.
This geometric condition is satisfied with $P=1$
for all $\xx,\yy\in\Omega_j$ if the region $\Omega_j$
is not much taller than $d$, and $\xx_0$ is placed near
the center of $\Omega_j$.
Thus, we will use $P=1$ in our numerical experiments.
However, if $\Omega_j$ is much taller than $d$
(high aspect ratio), then $P$ needs to be increased to guarantee
uniform convergence.
\end{remark}
   
Using this idea in the scattering setting
we represent the scattered field in the first layer as,
\begin{equation}
\label{repre1}
\us_1(\mathbf{x}) = \sum_{l= -P}^{P}\al^l\mathbf{S}^{k_1}_{\Gamma_1^l}\sigma_1+\sum_{l= -P}^{P}\al^l \mathbf{D}^{k_1}_{\Gamma_1^l}\mu_1 + \sum_{n\in\mathbb{Z}}a^{(1)}_n J_n(k_1\|\mathbf{x}-\mathbf{x}_1\|)e^{i n \theta_{\xx-\xx_1}}
~,
\qquad \mathbf{x}\in \Omega_1 
\end{equation}
where $\Gamma_1^l$ is the $l$th periodic translation of $\Gamma^0_1\in\Omega$, i.e. $\cup_{l\in\mathbb{Z}}\Gamma_1^l  = \Gamma_1$, $\sigma_1$ and $\mu_1$ are unknown periodic density functions defined on the interface $\Gamma_1$,
and $\xx_1\in\Omega_1$ is choice of origin.
The coefficients $\{a_n^{(1)}\}$ are now unknown and need to be determined by the boundary conditions. $\mathbf{S}$ and $\mathbf{D}$ are the usual single- and double-layer potentials \cite{Cot2}, which we may define living on a
general interface $\GG$ at wavenumber $k$ by,
\begin{align}
(\mathbf{S}^{k}_\GG\sigma)(\xx) & = \int_\GG G_{k}(\mathbf{x},\mathbf{y})\sigma(\mathbf{y}) ds_{\mathbf{y}}
\label{slp}\\
(\mathbf{D}^{k}_\GG\mu)(\xx) & = \int_\GG \frac{\partial G_{k}(\mathbf{x},\mathbf{y})}{\partial n(\mathbf{y})}\mu(\mathbf{y}) ds_{\mathbf{y}}
\label{dlp}
\end{align}
These representations satisfy the Helmholtz equation at wavenumber $k$
in $\mathbb{R}^2 \backslash \GG$;
thus the representation \eqref{repre1} satisfies the relevant Helmholtz equation in $\Omega_1$.
When restricted to target points on $\GG$ these
give the boundary integral operators $S^{k}_{\GG,\GG}$
(which is weakly singular), and $D^{k}_{\GG,\GG}$
(which is continous for $\GG$ smooth,
and is to be interpreted in the principal value sense).
Here and in what follows, the notation $S_{\Gamma_i,\Gamma_j}$ means the operator from source curve $\Gamma_j$ to target curve $\Gamma_i$.
We will also need the
operators corresponding to the target normal derivatives on $\GG$,
\begin{equation}
(N^{k}_{\GG,\GG}\sigma)(\xx) = \int_\GG \frac{\partial G_{k}(\mathbf{x},\mathbf{y})}{\partial n(\mathbf{x})}\sigma(\mathbf{y})ds_{\mathbf{y}}~, 
\qquad
(T^{k}_{\GG,\GG} \mu)(\xx) = \int_\GG \frac{\partial^2 G_{k}(\mathbf{x},\mathbf{y})}{\partial n(\mathbf{x}) \partial n(\mathbf{y})}\mu(\mathbf{y})ds_{\mathbf{y}}~,
\quad \xx\in\GG~.
\end{equation}
The operator $T^{k}_{\GG,\GG}$ is \textit{hypersingular} and defined in the Hadamard finite part sense.
The books \cite{Cot2,Cot} give further details.

We also need the {\em jump relations} that relate the limiting
values of \eqref{slp}--\eqref{dlp} to the actions of the above
boundary integral operators.
Let $u^{\pm}(\xx) := \lim_{h\to 0^+}u(\xx + h n(\xx))$
and $u_n^{\pm}(\xx) := \lim_{h\to 0^+} n(\xx) \cdot \nabla u(\xx + h n(\xx))$
be the limiting values and normal derivatives approaching $\xx\in\GG$
from the positive ($+$) or negative ($-$) side.
Then for all continuous densities $\sigma$ and $\mu$,
\bea
(\mathbf{S}^{k}_\GG\sigma)^\pm &=& S_{\GG,\GG}^{k} \sigma
\label{jr1}
\\
(\mathbf{S}^{k}_\GG\sigma)_n^\pm &=& (-\half + N_{\GG,\GG}^{k}) \sigma
\\
(\mathbf{D}^{k}_\GG\mu)^\pm &=& (\half + D_{\GG,\GG}^{k}) \mu
\\
(\mathbf{D}^{k}_\GG\mu)_n^\pm &=& T_{\GG,\GG}^{k} \mu
\label{jr4}
\eea
Thus the single-layer potential is continuous for all $\xx$,
whereas the double-layer is generally discontinuous across $\GG$.

Turning to the third layer $\Omega_3$, we similarly represent the scattered field $\us$ using layer potentials on $\Gamma_2$,
\begin{equation}
\label{repre2}
\us_3(\mathbf{x}) = \sum_{l= -P}^{P}\al^l\mathbf{S}^{k_3}_{\Gamma_2^l}\sigma_2+\sum_{l= -P}^{P}\al^l \mathbf{D}^{k_3}_{\Gamma_2^l}\mu_2 + \sum_{n\in\mathbb{Z}}a^{(3)}_n J_n(k_3\|\mathbf{x}-\mathbf{x}_3\|)e^{i n \theta_{\xx-\xx_3}}, \qquad \mathbf{x}\in\Omega_3
\end{equation}
where $\xx_3\in\Omega_3$.
The scattered field in the second layer has contribution from both $\Gamma_1$ and $\Gamma_2$, thus
\begin{align}
\label{repre3}
\us_2(\mathbf{x}) = &\sum_{l= -P}^{P}\al^l\mathbf{S}^{k_2}_{\Gamma_1^l}\sigma_1+\sum_{l= -P}^{P}\al^l \mathbf{D}^{k_2}_{\Gamma_1^l}\mu_1 \notag\\
& + \sum_{l= -P}^{P}\al^l\mathbf{S}^{k_2}_{\Gamma_2^l}\sigma_2+\sum_{l= -P}^{P}\al^l\mathbf{D}^{k_2}_{\Gamma_2^l}\mu_2 
+\sum_{n\in\mathbb{Z}}a^{(2)}_n J_n(k_2\|\mathbf{x}-\mathbf{x}_2\|)e^{i n \theta_{\xx-\xx_2}}, \quad
\mathbf{x} \in \Omega_2
\end{align} 
where $\xx_2\in\Omega_2$.

To determine the unknown densities $\sigma_1$, $\sigma_2$, $\mu_1$, $\mu_2$ and the coefficients $a^{(1)}_n$, $a^{(2)}_n$ and $a^{(3)}_n$, we enforce the following boundary conditions according to \eqref{cont} and \eqref{qp}--\eqref{rad2}.
\begin{itemize}
\item On $\Gamma^0_1$ and $\Gamma^0_2$, the continuity condition \eqref{cont} is imposed, giving,
\begin{align}
\label{bound1}
\begin{cases}
(\us_1-\us_2)|_{\Gamma_1^0} = -\ui|_{\Gamma_1^0} \\
\bigg(\frac{\partial \us_1}{\partial n}-\frac{\partial \us_2}{\partial n}\bigg)\bigg|_{\Gamma_1^0} = -\frac{\partial \ui}{\partial n}\bigg|_{\Gamma_1^0} \\
(\us_2-\us_3)|_{\Gamma_2^0}  = 0 \\
\bigg(\frac{\partial \us_2}{\partial n}- \frac {\partial \us_3}{\partial n}\bigg) \bigg|_{\Gamma_2^0}  = 0
\end{cases}
\end{align}
Note that since the representation of $\us$ in the three layers
involves layer potentials, the
limits must be taken from the appropriate side of $\Gamma_1$ and $\Gamma_2$
using jump relations \eqref{jr1}--\eqref{jr4}.

\item On $L_j$ and $R_j$, where $j = 1,2,3$, the quasi-periodicity
condition \eqref{quasi} is imposed:
\begin{align}
\label{bound2}
\begin{cases}
\al\us_j|_{L_j}- \us_j|_{R_j} = 0 \\
\al \frac{\partial \us_j}{\partial n }\bigg|_{L_j}- \frac{\partial \us_j}{\partial n }\bigg|_{R_j} = 0
\end{cases}
\end{align}
The left hand sides (phased differences)
is sometimes known as the {\em discrepancy} \cite{BG2011}.

\item On the parts of the artificial boundaries
$\Gamma_u$ and $\Gamma_d$ lying in the unit cell $\Omega$ (denote this part of the boundary by $\Gamma^0_u$ and $\Gamma^0_d$),
the radiation conditions \eqref{rad1} and \eqref{rad2} are imposed
for values and normal derivatives:
\begin{align}
\label{bound3}
\begin{cases}
(\us_1 - \sum_{n\in \mathbb{Z}}c_ne^{i\kappa_nx})|_{\Gamma^0_u}  = 0 \\
\bigg(\frac{\partial \us_1}{\partial n } - \sum_{n\in \mathbb{Z}}i c_n k_{1,n} e^{i\kappa_nx}\bigg)\bigg|_{\Gamma^0_u}  = 0 \\
(\us_3 - \sum_{n\in \mathbb{Z}}d_ne^{i\kappa_nx})|_{\Gamma^0_d}  = 0 \\
\bigg(\frac{\partial\us_1}{\partial n } + \sum_{n\in \mathbb{Z}}i d_n k_{3,n} e^{i\kappa_nx}\bigg)\bigg|_{\Gamma^0_d}  = 0
\end{cases}
\end{align}
\end{itemize}

\begin{figure}   
\centering
\includegraphics[height=52mm]{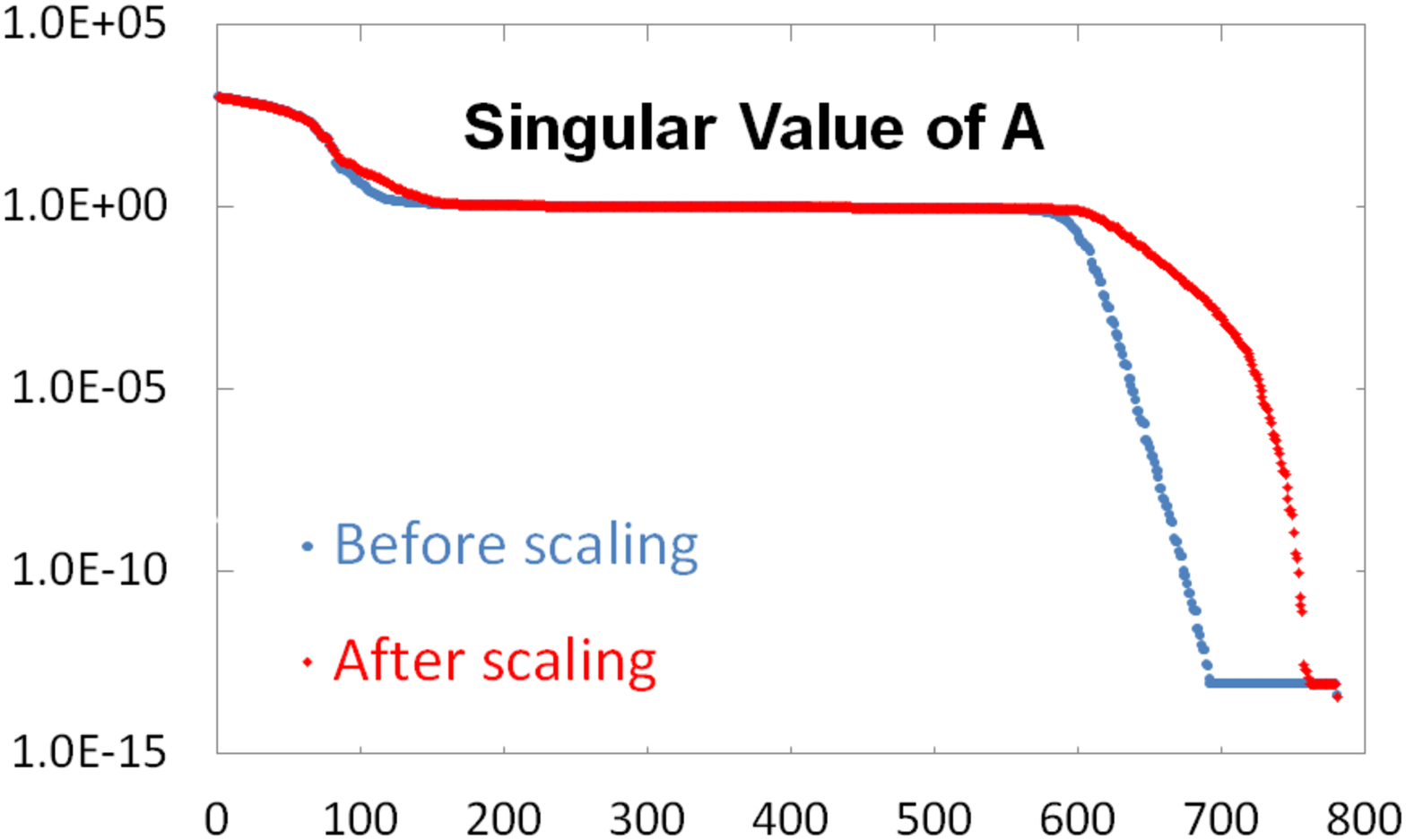}
\includegraphics[height=52mm]{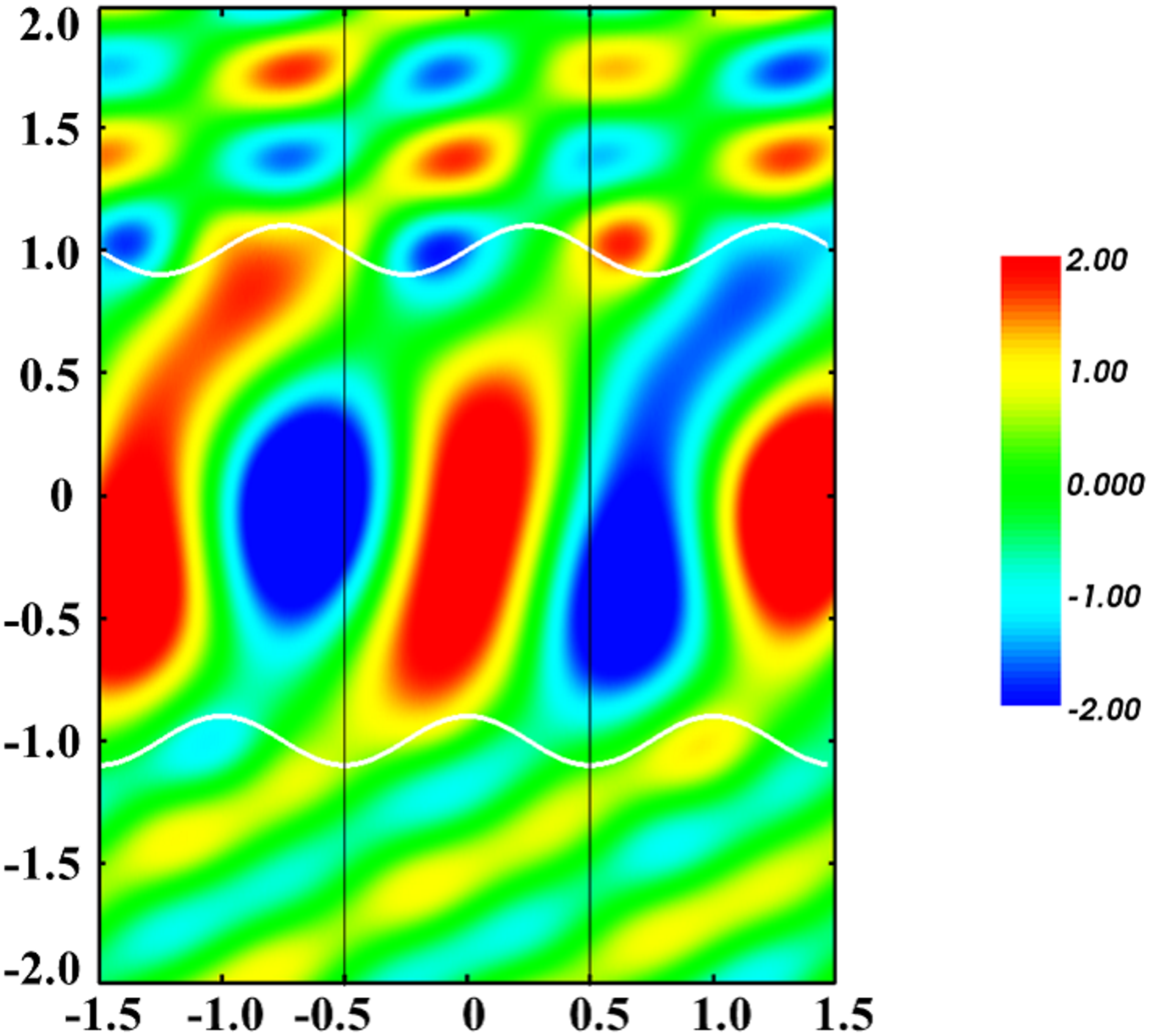}
\caption{Scattering from a periodic structure without inclusions, with $k_1=10$, $k_2=5$, $k_3=10$ and $d=1$. The interfaces are $\Gamma_1$ given by the graph $y=1+0.1\sin(2\pi x)$, and $\Gamma_2$ the graph $y=-1+0.1\cos(2\pi x)$. (a) The singular value spectrum $\sigma_j$ of matrix $\mathcal{A}$, vs index $j$:
original matrix (blue dots), vs after rescaling the columns of the blocks in $\mathcal{A}$ corresponding to $J$ expansions (red dots). (b) Real part of the scattered field computed to 13-digit accuracy.}
\label{singularvalue}
\end{figure}  

Substituting the representations \eqref{repre1}, \eqref{repre2} and
\eqref{repre3} into conditions \eqref{bound1}--\eqref{bound3}, one reaches
a system of coupled integral and functional equations
that can only be solved numerically, which means the interfaces must be discretized, and the infinite series truncated. 
The $J$ expansions are truncated up to order $Q$,
i.e.\ they retain $2Q+1$ terms. The centers of the $J$ expansion $\mathbf{x}_1$, $\mathbf{x}_2$ and $\mathbf{x}_3$ should be located roughly in the centers of the domains $\Omega_1$, $\Omega_2$ and $\Omega_3$.
We discretize the four interfaces $\Gamma_1$, $\Gamma_2$, $\Gamma_u$ and $\Gamma_d$ using equally spaced points ($\Gamma_1$ and $\Gamma_2$ are discretized through equally-spaced nodes in
their parametrizations). The left and right boundary $L_j$ and $R_j$ are discretized by Gauss--Legendre nodes.
The singular integrals involved in the layer potentials are discretized via
the Nystr\"om method \cite{Kress2010}, with
16th-order Alpert quadrature corrections \cite{Alpert1999}.
A phase correction is applied to account for the Bloch phase factors
when the parameter ``wraps'' around the end of the open curves
$\Gamma_1^0$ or $\Gamma_2^0$.
The standard application of quadrature rules, and the use of the Alpert
scheme, is described in \cite[Sec.~2.5]{mlqp} (see the smooth case only).

In the end, stacking as block rows the equations \eqref{bound1}, \eqref{bound2} and
\eqref{bound3}, the linear system for the discretized layered periodic structure without particle inclusions (the ``empty'' structure) takes the form,
\begin{equation}
\label{matA}
\left[\begin{array}{ccccccc}
A_{11} & A_{12}& A_{1m}& 0&0 &0 &0\\
A_{21} & A_{22}& A_{2m}& 0&0 &0 &0\\
A_{w1} & A_{w2}& A_{wm}& 0&0 &0 &0\\
A_{wu1}& 0&0 &A_{wuu} & 0&0 &0\\
0&A_{wd2} &0 &0 &0 &A_{wdd} &0\\
A_{u1} & 0 & 0 & A_{uu} & A_{ur} & 0 & 0\\
0&A_{d2} &0 &0 &0 &A_{dd} &A_{dr}\\
\end{array}
\right]
\left[\begin{array}{c}
\nu_1\\\nu_2\\\mbf{a}^{(2)}\\\mbf{a}^{(1)}\\\mbf{a}^{(3)}\\\mbf{c}\\\mbf{d}
\end{array} \right]
\; = \;
\left[\begin{array}{c}v^i\\0\\0\\0\\0\\0\\0
\end{array} \right] 
\end{equation}
Here the right-hand side vector has the form
\begin{equation}
f := [v^i,0, 0,0,0,0,0]^T \mbox{ with } v^i := [-\ui|_{\Gamma_1^0},\;-\partial \ui/\partial n |_{\Gamma_1^0}].
\label{Arhs}
\end{equation}
The unknown coefficient vector, which we will call $\valpha$,
stacks the discretizations of the paired
densities $\nu_j := [\mu_j;\sigma_j]$ for interfaces $j=1,2$,
the coefficient vectors $\mbf{a}^{(j)}$ for the $J$-expansions in layers
$j=1,2,3$,
and the Rayleigh--Bloch coefficient vectors $\mbf{c}$ and $\mbf{d}$
from \eqref{rad1}--\eqref{rad2}.
We may then summarize \eqref{matA} by
$$
\mathcal{A} \valpha = f
~.
$$


We now describe the matrix blocks in $\mathcal A$
(for readability we do not give formulae
for every single block, trusting that their construction is unambiguous from
the above; for more detail in a related scheme see \cite{mlqp}).
Each block maps unknowns to values and normal derivatives at target nodes,
or their phased differences between left and right walls.
Matrix entries involve either free-space Green's functions
between source and target nodes, or $J$ or Rayleigh--Bloch
expansions at target nodes.
\begin{itemize}
\item $A_{11}$ and $A_{22}$: Nystr\"om self-interaction matrices
for $\Gamma_1^0$ and $\Gamma_2^0$ respectively,
including the phased summation over $2P+1$ source near neighbors.
For instance,
\be
A_{11} = \left[\begin{array}{cc} I + \sum_{l=-P}^P \al^l (D^{k_1}_{\Gamma_1^0,\Gamma_1^l}-D^{k_2}_{\Gamma_1^0,\Gamma_1^l})
&\sum_{l=-P}^P \al^l (T^{k_1}_{\Gamma_1^0,\Gamma_1^l}-T^{k_2}_{\Gamma_1^0,\Gamma_1^l}) \\
\sum_{l=-P}^P \al^l (S^{k_1}_{\Gamma_1^0,\Gamma_1^l}-S^{k_2}_{\Gamma_1^0,\Gamma_1^l}) &
-I + \sum_{l=-P}^P \al^l (N^{k_1}_{\Gamma_1^0,\Gamma_1^l}-N^{k_2}_{\Gamma_1^0,\Gamma_1^l})
\end{array}\right]
\ee
(where here and below, the Nystr\"om discretizations of the various operators are implied),
maps $[\mu_1;\sigma_1]$ at the $N$
quadrature nodes on $\Gamma_1^0$ to the field jumps
$[\us_1-\us_2 ; \partial \us_1/\partial n - \partial \us_2/\partial n]$
at these same $N$ nodes. $I$ indicates the $N$-by-$N$ identity matrix.
All terms (after subtractions) are compact apart from the identities,
meaning that this subsystem is of Fredholm second kind.
Recall that the identities originate in the jump relations
\eqref{jr1}--\eqref{jr4}.
This BIE scheme for dielectric interfaces is due in the electromagnetic
case to M\"uller \cite{muller}, and the acoustic case to Kress--Roach
\cite{KR78} and Rokhlin \cite{rok2}.

\item $A_{21}$, $A_{12}$:
interaction of ($2P+1$ summed)
source densities on $\Gamma_1$ with value and normal
derivatives on $\Gamma_2$, and vice versa, at wavenumber $k_2$.

\item
$A_{1m}$, $A_{2m}$: values and normal derivatives
of the middle-layer $J$-expansion on $\Gamma_1^0$, $\Gamma_2^0$.

\item
$A_{w1}$, $A_{w2}$, and $A_{wm}$: phased differences of
values and normal derivatives between the left and right walls
$L_2$ and $R_2$ in \eqref{bound2} (case $j=2$),
due to the layer potentials on $\Gamma_1$, $\Gamma_2$, and the middle-layer
$J$-expansion respectively.
$A_{wu1}$, $A_{wuu}$:
phased differences for upper walls (case $j=1$) due to 
layer potentials on $\Gamma_1$ and the upper $J$-expansion.
$A_{wd2}$, $A_{wdd}$:
phased differences for lower walls (case $j=3$) due to 
layer potentials on $\Gamma_2$ and the lower $J$-expansion.

\item $A_{u1}$, $A_{uu}$ ($A_{d2}$, $A_{dd}$):
values and normal derivatives of the
layer potentials on $\Gamma_1$ ($\Gamma_2$), and the upper (lower) $J$-expansions, evaluated on the upper boundary $\Gamma^0_u$ (lower boundary $\Gamma^0_d$).

\item $A_{ur}$, $A_{dr}$: values and normal derivatives of the Rayleigh--Bloch expansions on the upper and lower boundaries $\Gamma^0_u$, $\Gamma^0_d$
respectively.
\end{itemize}

\begin{remark}
The lower-right $5\times 5$ block of $\mathcal A$ in \eqref{matA} involves
only the effect of the auxiliary periodizing degrees of freedom
($J$ and Rayleigh--Bloch expansions) on the auxiliary matching
conditions (discrepancies and radiation conditions).
In related work this block is given the symbol $Q$ \cite{BG2011,mlqp}.
\end{remark}

The matrix $\mathcal{A}$ is generally rectangular, depending on the specific
numbers of discretization nodes.
Although the upper-left $2\times 2$ block of $\mathcal A$ is
a square system coming from a Fredholm second kind system of BIEs,
the other blocks involve $J$-expansions evaluated on interfaces
and walls, which make $\mathcal A$ as a whole exponentially ill-conditioned.
Figure \ref{singularvalue}(a) shows the singular values of $\mathcal{A}$: there are many singular values clustered around $10^{-14}$, although the situation is alleviated a little by rescaling the columns that correspond to the $J$-expansions (also shown).
However, such ill-conditioning is not an obstacle as long as the system
is consistent:
since $\mathcal A$ is not too large (of typical size $10^3$ for
$k_j$ corresponding to up to several wavelengths across one period $d$),
we may use direct dense linear algebra for a small-norm least-squares solution.
We use the {\tt mldivide} command in MATLAB.
Figure \ref{singularvalue}(b) shows the resulting real part of the scattered field with a total flux error (whose definition will be given in section \ref{Numerical}) of $10^{-13}$.

\section{Multi-particle scattering}
\label{scat_multi}

Multi-particle scattering in free space has been discussed extensively in \cite{GG2013,Lai:14,ALM02}, with applications including climatology, remote sensing, and design of composite materials. Recently, we have developed a fast solver for finding the field scattered from a large number of particles located in a layered medium, by combining the Sommerfeld integral and multiple scattering theory. Here we briefly review the method introduced in \cite{Lai:14}, and then combine
it with the periodic grating. See \cite{Lai:14,GG2013} for more details.
Note that our approach can be seen as a simple
version of a reduced basis method \cite{chenRBM}.
 
\subsection{Scattering matrix of a single particle}
Consider for now a dielectric inclusion with wavenumber $k_p = \omega\sqrt{\epsilon_p \mu}$ surrounded by uniform dielectric with $k_2 = \omega\sqrt{\epsilon_2 \mu}$.  When the inclusion is a disk of radius $R$ centered at the origin, 
it is well known that the solution can be represented using separation of variables, with
\begin{equation}\label{int}
u(r,\theta) = \sum_{n=-\infty}^{\infty} \alpha_n J_n(k_p r)e^{i n\theta}
\end{equation}
in the interior, and scattered wave
\begin{equation}\label{out}
\us(r,\theta) = \sum_{n=-\infty}^{\infty} \beta_n H^{(1)}_n(k_2 r)e^{i n\theta}
\end{equation}
in the exterior. Here, $(r,\theta)$ are the polar coordinates of a point in the plane, and $H^{(1)}_n(r)$ is the Hankel function of the first kind of order $n$.

\begin{definition}
The mapping between the incoming coefficients $\{\alpha_n\}$ and outgoing coefficients $\{\beta_l\}$ is referred as the {\em scattering matrix}.
It will be denoted by $S$, with matrix elements $s_{ln}$.
\end{definition}

For a single disk, the scattering matrix $S$ is diagonal and is easily found
analytically \cite{Cot2}; this is not true for an arbitrary inclusion shape.
We instead seek a solution via BIE, using the M\"uller--Kress--Roach--Rokhlin scheme from the previous section.
Suppose that the inclusion $\Omega_p$ has boundary $\partial \Omega_p$
and is enclosed by a disk $D$ centered at the origin.
Given the incident wave $\ui$ and the boundary conditions \eqref{cont}, the exterior scattered field $\us$ and the field $u$ within $\Omega_p$ have the following representations \cite{Cot2}:
\begin{align}
u^s&= \mathbf{S}^{k_2}\sigma+\mathbf{D}^{k_2}\mu, \qquad \mbox{ for } \mathbf{x}\in \Omega_p^c,\label{rep1}\\
u& = \mathbf{S}^{k_p}\sigma+\mathbf{D}^{k_p}\mu, \qquad \mbox{ for }\mathbf{x}\in \Omega_p. \label{rep2} 
\end{align} 
Here $\sigma$ and $\mu$ are unknown single- and double-layer densities on
$\partial \Omega_p$, and, in this section without ambiguity we drop the
subscripts $\partial \Omega_p$.
Enforcing the interface conditions \eqref{cont} and taking appropriate limits
using jump relations \eqref{jr1}--\eqref{jr4}
yields a system of Fredholm integral equations of the second kind:
\begin{align}
\mu+ (S^{k_2}-S^{k_p})\sigma
+(D^{k_2}-D^{k_p})\mu & = -\ui, \label{intg1}\\
-\sigma+(N^{k_2}-N^{k_p})\sigma
+(T^{k_2}-T^{k_p})\mu & = -\frac{\partial \ui}{\partial n}.
\label{intg2}
\end{align}

%



Let $\sigma_n$ and $\mu_n$ denote the solution to
\eqref{intg1}--\eqref{intg2} for
$\ui(r,\theta) = J_n(kr)e^{i n \theta}$.
We may then precompute the scattering matrix elements $s_{ln}$
as the multipole expansion coefficients (truncated up to $2p+1$ terms)
\begin{equation}
\us(r,\theta) \approx \sum_{l=-p}^{p} s_{ln} H^{(1)}_l(k_2 r)e^{i l\theta} \, ,
\end{equation}
from the densities $\sigma_n$ and $\mu_n$ via Graf's addition theorem,
giving the standard formula \cite{rok90,wideband}
\begin{equation}
s_{ln}
= \int_{\partial \Omega_p} 
J_l(k_2\|\mathbf{y}\|) e^{- i l \theta_{\mathbf{y}}} \, 
\sigma_n(\mathbf{y}) + \, 
n(\yy) \cdot \nabla [J_l(k_2\|\mathbf{y}\|) e^{- i l \theta_{\mathbf{y}}}]
\mu_n(\mathbf{y})
\;  ds_{\mathbf{y}} \, .
\end{equation}

\begin{figure}[tbp]  
\centering
\includegraphics[width=70mm]{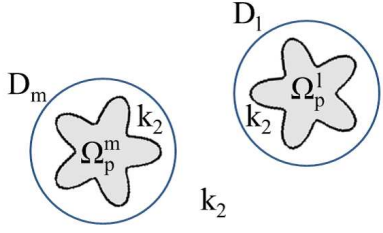}
\caption{Two inclusions and their enclosing disks. 
The scattering matrix $S_j$ for each inclusion $\Omega^j_p$ 
with wavenumber $k_p$ is defined as the map from an incoming
field on $D_j$ to the corresponding outgoing field. 
}
\label{scatdisk}
\end{figure}

\subsection{Multiple inclusions}

Suppose now that we have $M$ inclusions $\Omega^1_p,\dots,\Omega^M_p$ that are identical up to translation and rotation, and are well separated in the sense that each inclusion $\Omega^m_p$ lies within a disk $D_m$ of radius $R$ such that the disks are not overlapping (see \fref{scatdisk}).
The incident wave for the $m$th inclusion may now be expanded as 
\be
\ui(\xx) \approx \sum_{n=-p}^{p} a^{(m)}_n J_n(k_2 r_m)e^{i n\theta_m}
\label{incident}
\ee
where $(r_m,\theta_m)$ are defined to be the polar coordinates
of $\xx$ relative to the center of $D_m$.
We will denote by $\valpha^m$ the set of
$2p+1$ incoming coefficients, and by $\vbeta^m$ 
the $2p+1$ outgoing coefficients, for the $m$th particle. 
Thus
\begin{equation}
\vbeta^m = S_p^{(m)} \valpha^m, \qquad \mbox{ for } m = 1,\dots,M.
\label{scat}
\end{equation}
where $S_p^{(m)}$ denotes the truncated $(2p+1) \times (2p+1)$ 
scattering matrix acting on the truncated expansion about the center of the
$m$th particle.
If the particles were only translated, we would have  $S_p^{(m)} =S_p$ for all $m$,
where $S_p$ is the truncated scattering matrix from the previous section
with elements $s_{ln}$.
We allow general rotations of particles;
rotation of the $m$th particle by angle $\phi$ introduces phase factors, so
$(S_p^{(m)})_{ln} = e^{i\phi(l-n)}s_{ln}$.

Multi-particle scattering has a key difference from single particle scattering, namely that the incoming field experienced by each particle consists of two parts: the (applied) incident field $\ui$,
and the contribution to the scattered field $\us$ from all of the {\em other} particles.
We denote by 
$T^{jm}$ the {\em translation operator}
(or M2L in FMM language \cite{rok90})
that maps the outgoing coefficient vector  $\vbeta^m:=\{\beta_n^m\}_{n=-p}^p$ 
from particle $m$ to their contribution to the local expansion coefficients $\valpha^j$ centered at particle $j$.
With this operator in place, the 
incoming coefficients $\valpha^m$ for the $m$th particle are
\begin{equation}\label{trans}
\valpha^m = \veca^m+\sum_{\substack{j = 1 \\ j\neq m}}^{M} T^{jm}\vbeta^j,
\end{equation}
where $\veca^m$ is the (truncated) local expansion \eqref{incident}
of the incident wave $\ui$ relative to particle $m$. 

Combining \eqref{scat} and \eqref{trans}, one can easily eliminate the 
incoming coefficients $\valpha^m$ to obtain the following linear system that 
only involves the outgoing coefficients: 
\begin{equation}\label{scattmatrix}
\left(
\mathcal{S}^{-1}-\mathcal{T}\right)
\left[\begin{array}{c}
\vbeta^1 \\
\vbeta^2 \\
\vdots \\
\vbeta^M
\end{array} 
\right] 
=
\left[\begin{array}{c}
\veca^1 \\
\veca^2 \\
\vdots \\
\veca^M
\end{array} 
\right],
\end{equation}
where 
\begin{equation*}
\mathcal{S} := \left[\begin{array}{cccc}
S_p^{(1)} & & &\\
& S_p^{(2)} & &\\
& & \ddots &\\
& & & S_p^{(M)}
\end{array}
\right],
\quad
\mathcal{T} := \left[\begin{array}{cccc}
0 & T^{21} &\cdots &T^{M1}\\
T^{12} & 0 &\cdots &T^{M2}\\
\vdots&\vdots & \ddots &\vdots\\
T^{1M}&T^{2M} &\cdots & 0
\end{array}
\right].
\end{equation*}

The system \eqref{scattmatrix} can be solved iteratively, using GMRES. Since each translation operator $T^{nm}$ is dense, a naive matrix-vector product requires $\bigO(M^2(2p+1)^2)$ operations, where $p$ is the order of the truncated expansion. The cost can be reduced to $\bigO(M(2p+1)^2)$ work by FMM acceleration, for which we refer the reader to \cite{rok90,wideband}. Furthermore, there exists an effective preconditioner for the system \eqref{scattmatrix}. Left-multiplying by the block diagonal matrix $\mathcal{S}$ results in the preconditioned system matrix  $I-\mathcal{ST}$. This significantly reduces the number of iterations.

The advantage of using the one-particle scattering matrices $S_p^{(m)}$ over
boundary integral equations is clear: the number of degrees of freedom per inclusion is only $2p+1$ rather than the number of nodes needed to discretize the domain boundaries $\partial \Omega^m_p$. For complicated inclusions, this permits a vast reduction in the number of degrees of freedom required, forming the basis for the so-called FMPS method \cite{GG2013}. Moreover, the block-diagonal preconditioned multiple scattering equations are much better conditioned than the BIE \eqref{intg1}--\eqref{intg2}, while FMM acceleration is particularly fast in this setting.

\begin{remark}
It is straightforward to extend the method to more than one type of particles as long as the assumption that the enclosed circles are well separated still holds. The additional cost is simply the bookkeeping for the different scattering matrices.
\end{remark}

\section{Multi-particle scattering in the periodic layered medium}
\label{Quasi_scheme}

We now combine the schemes of the previous two sections.
The field in the middle layer is the periodic layered contribution $\us_2$ from
\eqref{repre3}, plus the scattered field from the $M$ inclusions
and, cruicially, their $2P$ neighboring near-field phased copies.
We need the notation $(r_m^l,\theta_m^l)$ for polar coordinates
relative to the origin of the $m$th particle translated by $(ld,0)$.
Then,
\begin{align}
\label{repre3_mod}
u_2(\mathbf{x}) = \us_2(\mathbf{x}) + \sum_{l=-P}^P \alpha^l
\sum_{j=1}^M \sum_{n=-p}^p 
\beta_n^m H^{(1)}_n(k_2r^l_m) e^{i n \theta^l_m}, \quad \mathbf{x} \in
\Omega_2\backslash \Omega_p ~.
\end{align}    

It only remains to set up the interactions between
the layered periodic structure and the inclusion structure.
We denote the translation matrix mapping layer densities and the middle-layer $J$ expansion to the incoming coefficients of all particles by $\mathcal{B}$, and the translation matrix mapping (phased summed) outgoing particle
multipole coefficients to data on the layers and walls by $\mathcal{C}$.
Adhering to the ordering of unknowns and conditions in \eqref{matA},
one can show that they have the forms (since the particles interact only with
the middle layer),
\bea
\label{matB}
\mathcal{B} = \left[\begin{array}{ccccccc}
B_{p1} & B_{p2}& B_{pm}& 0&0 &0 &0
\end{array}
\right],
\\
\label{matC}
 \mathcal{C} = \left[\begin{array}{ccccccc}
 C_{1p} & C_{2p} & C_{wp} &
 0 &  0 &  0 & 0
\end{array}
\right]^T,
\eea
where here $[\cdot]^T$ denotes the blockwise transpose.
$B_{p1}$, $B_{p2}$ and $B_{pm}$ map layer densities on $\Gamma_1$, $\Gamma_2$ and the $J$ expansion to the incoming coefficients of the particles.
$C_{1p}$, $C_{2p}$, and $C_{wp}$ map the outgoing coefficients from particles to the values and normal derivatives on interfaces $\Gamma_1$, $\Gamma_2$ and the
discrepancy from $L_2$ to $R_2$.

It is easy to construct the elements in $C_{1p}$, $C_{2p}$, and $C_{wp}$
by direct evaluation of multipole expansions.
Note that $C_{wp}$ involves {\em cancellations} which mean that
only particle-wall interactions over distances greater than $d$ survive,
as discussed in Remark~\ref{cancel}.

To obtain the elements in $B_{p1}$, $B_{p2}$ and $B_{pm}$, one again uses
Graf's addition theorem.
In particular, the translation submatrix $B_{pm}$ that maps the coefficients from one $J$ expansion to another $J$ expansion (the local-to-local or L2L operator in the FMM) is constructed through the following lemma.

\begin{lemma}[\cite{rok90}]
Let disk $m$ be centered at $\mathbf{x}_m$
and let disk $l$ be centered at $\mathbf{x}_l$.
Then the local expansion 
\begin{equation}
\sum_{n=-\infty}^{\infty} \gamma_n^m J_n(k_2 r_m)e^{i n\theta_m}
\end{equation}
induces a field on disk $l$ of the form
\begin{equation}
u = \sum_{n'=-\infty}^{\infty} \alpha^l_{n'} J_{n'}(k_2 r_l)e^{i n'(\theta_l-\pi)}
\end{equation}
where
\[
\alpha^l_{n'} = \sum_{n=-\infty}^{\infty} 
e^{i(n-n')\theta_{\mathbf{x}_m -\mathbf{x}_l}} \gamma^m_{n-n'}
				J_n(k_2 \| \mathbf{x_m} - \mathbf{x_l} \|).
\]
\end{lemma}

Let us denote by $\vbeta:=\{\vbeta^m\}_{m=1}^M$ the multipole coefficients for all $M$ particles in $\Omega_2$. Combining the matrices \eqref{matA}, \eqref{scattmatrix}, \eqref{matB} and \eqref{matC}, we obtain the final system:
\begin{equation}
\label{finalsys}
\left[\begin{array}{cc}
\mathcal{A} & \mathcal{C} \\
\mathcal{B} & \mathcal{D} 
\end{array}
\right]
\left[\begin{array}{c}
\valpha \\
\vbeta
\end{array}
\right]
=
\left[\begin{array}{c}
f \\
0
\end{array}
\right],
\end{equation}
where $\mathcal{D} = \mathcal{S}^{-1}-\mathcal{T} $,
as in \eqref{scattmatrix}, and $f$ is the right-hand side vector \eqref{Arhs}.

Since \eqref{finalsys} is a rectangular ill-conditioned matrix, we
cannot easily solve this whole system iteratively. 
However, we now present a Schur complement scheme to
generate a smaller, well-conditioned square linear system for which
an iterative solution is efficient.
Since $\mathcal A$ has size of order $10^3$ in both dimensions, it is much smaller than $\mathcal D$,
and we can eliminate the unknowns $\valpha$ via
$\valpha = \pA (f - C\vbeta)$,
where $\pA$ is the pseudoinverse of $\mathcal A$.
We precompute $\mathcal A = U\Sigma V^\ast$, the singular value decomposition
of $\mathcal A$, where the singular values are $\sigma_j$.
Then to apply $\pA$ to an arbitrary vector $g$ we use matrix-vector multiplies,
\be
\pA g = V \Sigma^{+} (U g)~,
\label{applypA}
\ee
where $\Sigma^{+}$ has diagonal elements $\min[1/\sigma_j, 1/\epsilon]$.
The regularization parameter is fixed at $\epsilon = 10^{-10}$;
its precise choice is not crucial, but empirically it is best chosen to be
roughly the desired solution accuracy.
\begin{remark} Forming the matrix $\pA$ then multiplying it against a vector $g$ is dangerous,
since it can cause large and unnecessary round-off error.
Rather, \eqref{applypA} should be used since it is backward stable
and hence introduces the minimum possible round-off error.
\end{remark}

Substitution for $\valpha$ into the full linear system gives the
Schur complement
\begin{equation}\label{schur}
(\mathcal{D}-\mathcal{B}\pA\mathcal{C})\vbeta = -\mathcal{B}\pA f~,
\end{equation}
a relatively well-conditioned square system involving only unknowns $\vbeta$.
The new system matrix has a physical interpretation:
it is the particle-particle interaction matrix using the layered-medium
quasi-periodic Green's function, where
$\mathcal{B}\pA\mathcal{C}$ is a low-rank update to $\mathcal D$.
We can use GMRES on this system,
using FMM acceleration to apply $\mathcal D$,
$\mathcal{B}$ and $\mathcal{C}$,
and the factorization
\eqref{applypA} to apply $\pA$.
The cost (for moderate frequencies) is then an optimal
$\bigO(M)$ per iteration.

As we mentioned in section 4, left-preconditioning via $\mathcal S$
can improve the conditioning of the system. We therefore end up solving the following,
\begin{equation}
(I - \mathcal{S}\mathcal{T}-\mathcal{S}\mathcal{B}\pA
\mathcal{C})\vbeta = -\mathcal{S}\mathcal{B}\pA f ~,
\end{equation}
with the same cost per iteration as \eqref{schur}
(the small dense multiplications by $S_p$  being cheap),
but fewer needed iterations.

\section{Numerical experiments}
\label{Numerical}

In this section, we demonstrate the performance of our algorithm with three examples. For simplicity, we use a single class of particles, parametrized by 
\begin{equation}\label{par1}
\left\{\begin{array}{rcl}
x &=& (a_1+a_2\cos(a_3 t))\cos(t),\\
y &=& (a_1+a_2\cos(a_3 t))\sin(t),
\end{array}\right.\mbox{ for } 0\leq t < 2\pi.
\end{equation}
Particles with more complicated boundaries do not introduce any essential difficulty in our scheme, except that the precomputation of the scattering matrix is a little more involved, particularly if corners are present \cite{BRS2010,Helsing2008}. However, regardless of complexity, the 
size of the needed scattering matrix $(2p+1)$
depends only on the particle size in wavelengths and the closeness of
nearby particles.
Given a fixed $a_1$, $a_2$ and $a_3$, multiple copies of the inclusion are randomly distributed in the central layer of the medium with random orientations. 
  
To provide an independent test of accuracy, we check flux conservation. If all wavenumbers $k_1$, $k_2$, $k_3$ and $k_p$ are real, the Rayleigh--Bloch coefficients satisfy the identity (eg see \cite{Linton2007165})
\begin{equation}
\sum_{k_{1,n}>0} k_{1,n}|c_i|^2 +\sum_{k_{3,n}>0}k_{3,n}|d_i|^2= k_1\sin\theta
\end{equation}
In other words, the outgoing energy flux must equal the incoming energy flux.
{\em Flux error} then refers to the size of the
difference between the left and right sides.

Throughout all the numerical examples, the period is $d=1$,
and the near-field summation is fixed at $P=1$ (i.e.\ three
terms in the near-field sum).
In the first two examples the interfaces $\Gamma_1$ and $\Gamma_2$ are given by
the graphs $y = 1 + 0.1\sin(2\pi x)$ and $y = -1 + 0.2\cos(2\pi x)$
respectively.

Assume there is a modest distance between the inclusions $\Omega_p$ and the interfaces $\Gamma_1$ and $\Gamma_2$, say, at least $0.5$ wavelengths in terms of the wavenumber $k_2$. Under this assumption, we may discretize $\Gamma_1$ and $\Gamma_2$ with $N=120$ nodes each, equally spaced in $x$,
which is sufficient to achieve $12$ digits of accuracy  with the periodized Alpert's quadrature for modest wavenumbers $k_1$, $k_2$ and $k_3$.
The two artificial boundaries $\Gamma_u$ and $\Gamma_d$ are discretized equally in $x$ with 50 points each.
We use $2Q+1 =73$ terms for the $J$ expansions in each layer, and $41$
Rayleigh--Bloch modes in each vertical direction. Based on these discretizations, $\mathcal{A}$ has size $1020\times 781$.

All computations are carried out using a 2.3GHz Intel Core i5 laptop, with 4GB RAM.

\subsection{Example 1: scattering from large numbers of inclusions }
\label{NumericalExperiment1}

\begin{figure}[tbp]  
\centering
\includegraphics[width=140mm]{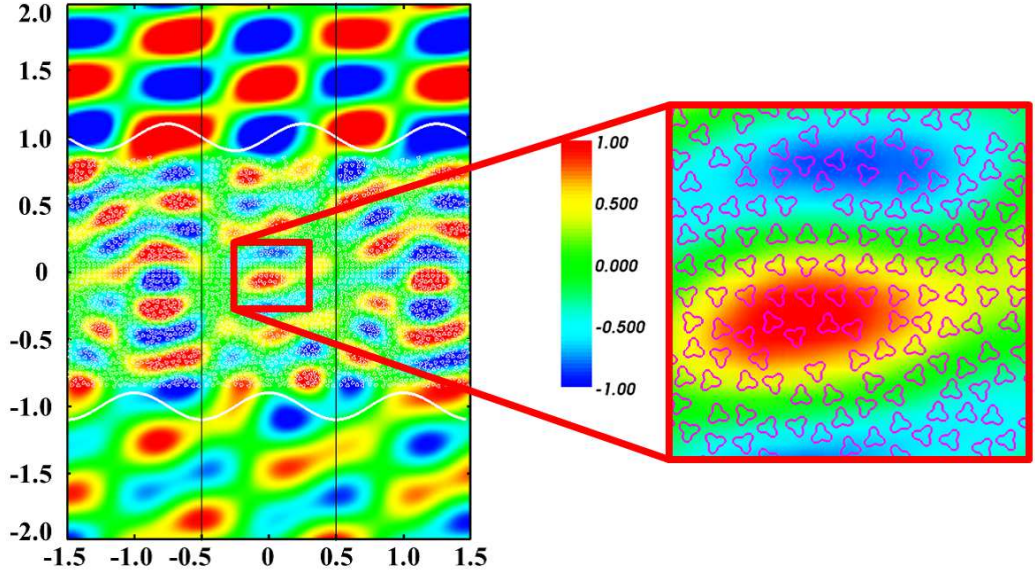}
\caption{Real part of the total field with $1000$ dielectric inclusions randomly distributed in a three-layered medium (see Example 1).
The wavenumber for each particle is $k_p=30$ and the wavenumbers for the three layers are $k_1=10$, $k_2=8$, $k_3=10$. 
The diameter of each particle is approximately $0.2$ wavelengths at the wavenumber $k_p$.}
\label{fig:NumericalResult1}
\end{figure}

\begin{figure}[tbp]   
\centering
(a)\raisebox{-1.3in}{\includegraphics[width=68mm]{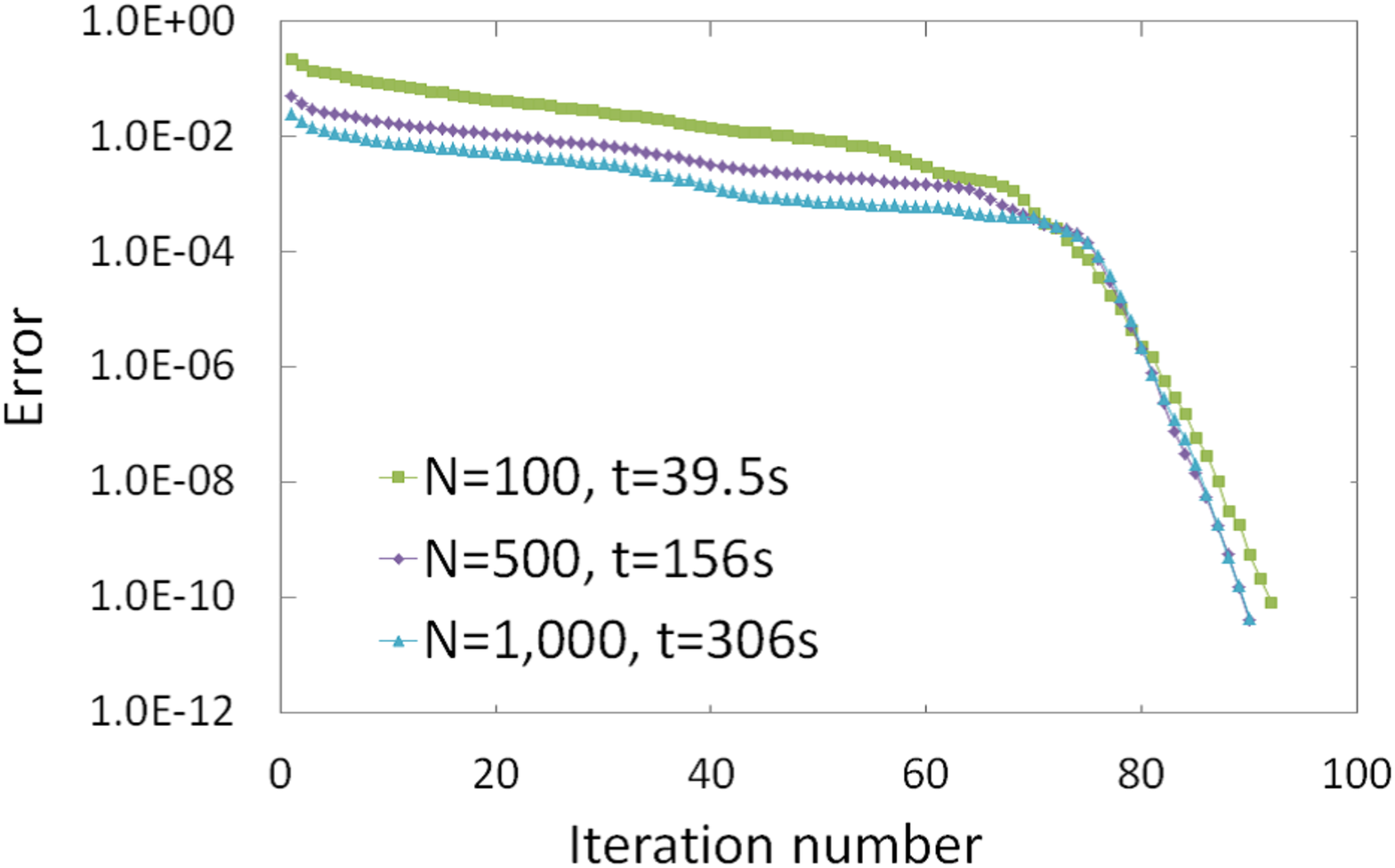}}
\quad
(b)\raisebox{-1.3in}{\includegraphics[width=68mm]{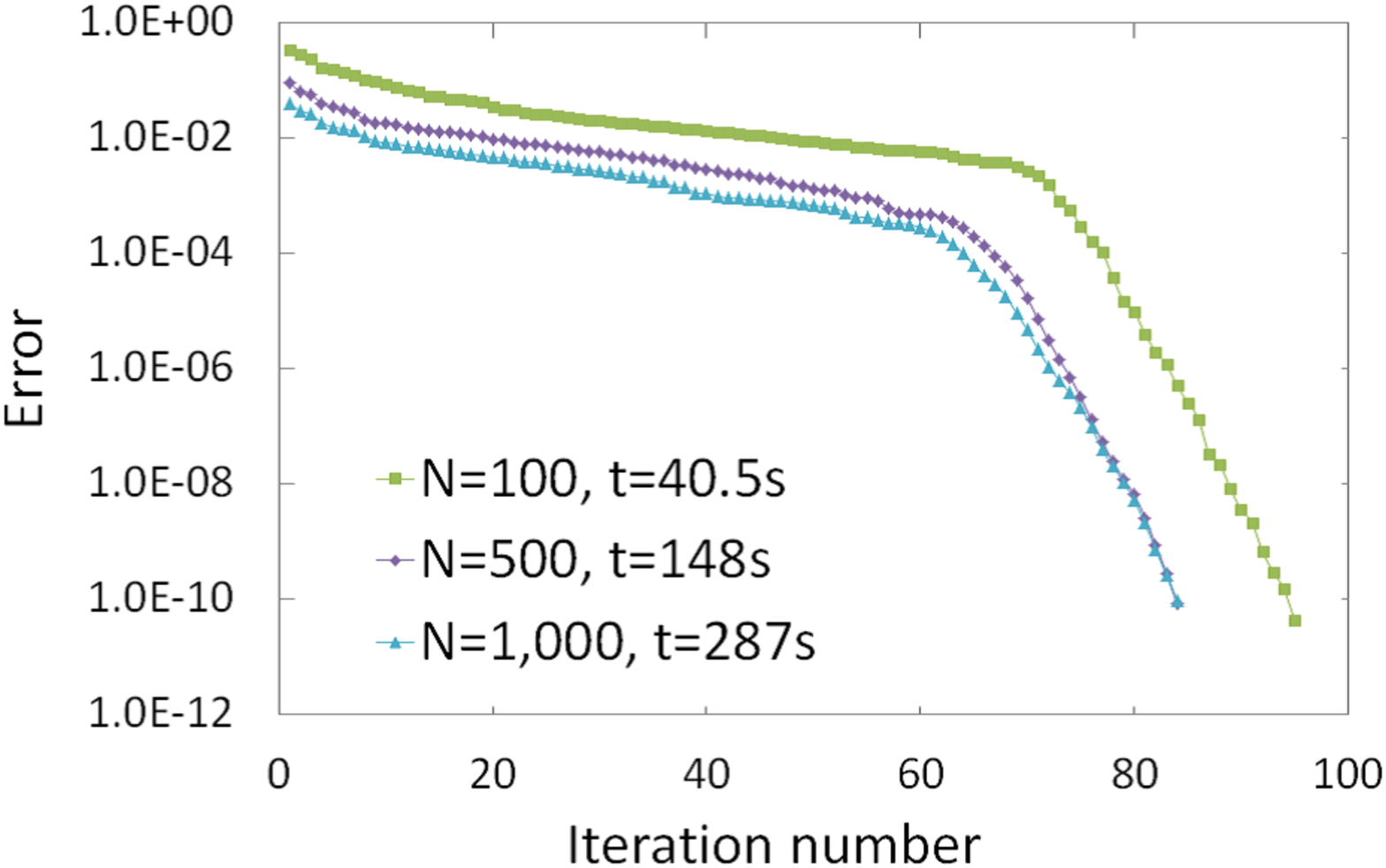}}
\caption{Convergence behavior of GMRES and the CPU time required for various numbers of periodized inclusions embedded in either (a) homogeneous medium or (b) a three-layered periodic medium, for Example 1.
For (a), we set $k_1=k_2=k_3=10$ and for (b), we set $k_1=10$, $k_2=8$, $k_3=10$.}
\label{fig:NumericalResult1-1}
\end{figure}

In our first example, we consider distributions of $M= 100,\ 500$ and $1000$ inclusions with wavenumber $k_p = 30$, by assuming the incident angle is away from any Wood's anomaly, and the boundaries of the inclusions do not touch the wall $L_2$ or $R_2$. (Both restrictions will be removed, in our second and third examples, respectively.)
The thickness of the central layer is fixed by the length of $L_2$ and $R_2$, given by $l=2.0$.
The size of each inclusion is therefore determined by the amount. Assume $a_3=3$ in Eq. \eqref{par1}.
We let $a_1=0.0309,\ a_2=0.0103$ for $M=100$, and $a_1=0.0154,\ a_2=0.00514$ for $M=500$, and $a_1=0.0111,\ a_2=0.0037$ for $M=1000$.
To obtain the scattering matrix $S_p$ with $p = 10$, we solve the integral Eqs. \eqref{intg1} and \eqref{intg2} by discretizing the boundary of the particle with $N=300$ equispaced points. We assume the wavenumbers of the layered medium are given by $k_1=10,\ k_2=8$ and $k_3=10$.
The incident angle is set to be $\theta=-\arccos(1.0-2\pi /10.0)+0.1\approx  -1.089976736488571$, which is not a Wood's anomaly for $k_1 = 10$, but is
quite close to one.
Results are presented in Figs. \ref{fig:NumericalResult1}, \ref{fig:NumericalResult1-1} and Table~\ref{tab:NumericalResult1}. 

Figure \ref{fig:NumericalResult1} shows the total field in the case $M=1000$.
Disordered propagation due to the random inclusions is apparent.
The total number of unknowns in $\vbeta$ is $21000$, although if
nodes on the particles were used it would be much higher.
It requires $287$s to achieve $9$ digits of accuracy. 
Figure \ref{fig:NumericalResult1-1} shows the convergence behavior of GMRES as the number of inclusions is increased, and the total CPU times. In Fig.~\ref{fig:NumericalResult1-1}(a), we also study the convergence rate when the background is homogeneous, by setting the material parameters to be the same for the three layers ($k_1=k_2=k_3=10$). No obvious discrepancy in terms of the number of iterations has been observed, which suggests that the multiple scattering is dominated by the inclusions.
Table \ref{tab:NumericalResult1} shows flux error in the cases above.
In all cases we exceed $9$ digits of accuracy.

\begin{table}[tbp] 
\begin{center}
\begin{tabular}{|c|c|c|}
\hline
Number of particles & flux error, homogeneous case & flux error, three-layer case\\
\hline
$100$ & $1.05$e-9 & $2.22$e-10 \\
$500$ & $4.03$e-10 & $2.35$e-10 \\
$1,000$ & $1.59$e-9 & $1.14$e-9 \\
\hline
\end{tabular}
\end{center}
\caption{Flux error for various numbers of inclusions embedded in either (a) homogeneous medium or (b) a three layered medium, for Example 1.}
\label{tab:NumericalResult1}
\end{table}

\subsection{Example 2: scattering from large numbers of inclusions at a Wood's anomaly}
 \label{NumericalExperiment2}

\begin{figure}[tbp] 
\begin{center}
\includegraphics[width=140mm]{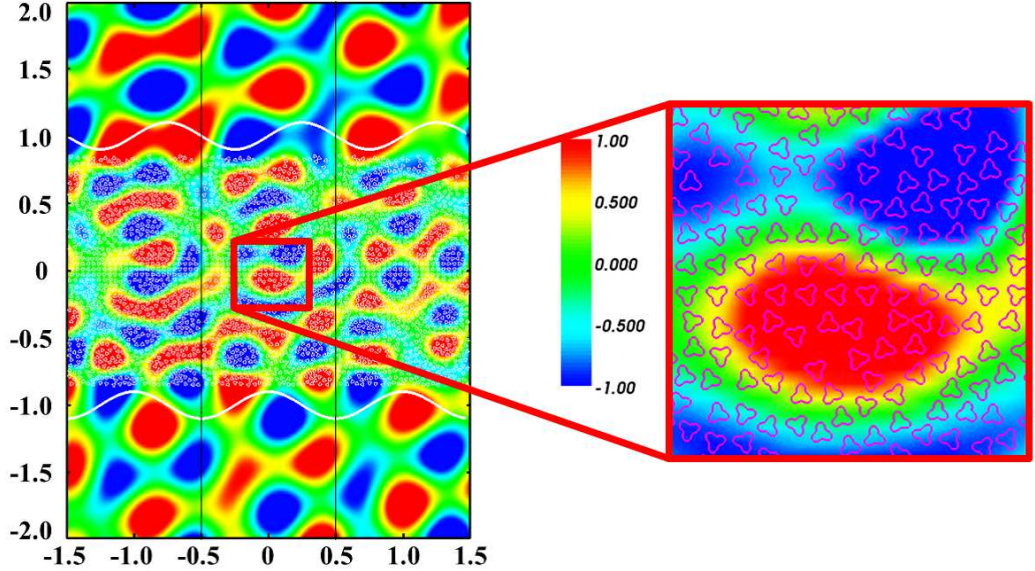}
\caption{Real part of the total field when $1000$ dielectric inclusions with $k_p=30$ are randomly embedded in a three-layered medium with $k_1=10$, $k_2=8$ and $k_3=10$, as in \fref{fig:NumericalResult1-1}.
This is Example 2: the angle of the incident wave is at a Wood's anomaly.
}
\label{fig:NumericalResult2}
\end{center}
\end{figure}

\begin{table}[tbp]
\begin{center}
\begin{tabular}{|c|c|c|c|c|}
\hline
Number of particles & number of iterations & GMRES error & CPU time (sec.) & Flux error \\
\hline
$100$ & $97$ & $7.61$e-11 & $41.2$ & $1.03$e-9 \\
$500$ & $91$ & $4.24$e-11 & $156$ & $3.21$e-10 \\
$1,000$ & $90$ & $5.83$e-11 & $307$ & $2.48$e-9 \\
\hline
\end{tabular}
\end{center}
\caption{Convergence behavior of GMRES, the CPU time required and flux error for various numbers of inclusions embedded in a three layered medium (see Example 2.)}
\label{tab:NumericalResult2}
\end{table}

In our second example, we consider the same scattering as above except the change of incident angle. In particular, we let $\theta=-\arccos(1.0-2\pi /10.0)\approx -1.189976736488571$, which is a Wood's anomaly for $k_1=10$.
Thus the classical quasi-periodic Green's function for the upper layer
does not exist.
However, our scheme is still able to obtain $10$ digits of accuracy as shown in Fig. \ref{fig:NumericalResult2} and Table \ref{tab:NumericalResult2}.
The flux error suggests at least $9$ digits of accuracy in all cases. The number of GMRES iterations and CPU times are almost the same as in Example 1.

\subsection{Example 3: scattering from inclusions intersecting with unit cell walls}
\label{hitwalls}

\begin{figure}[tbp]  
\begin{center}
\includegraphics[width=140mm]{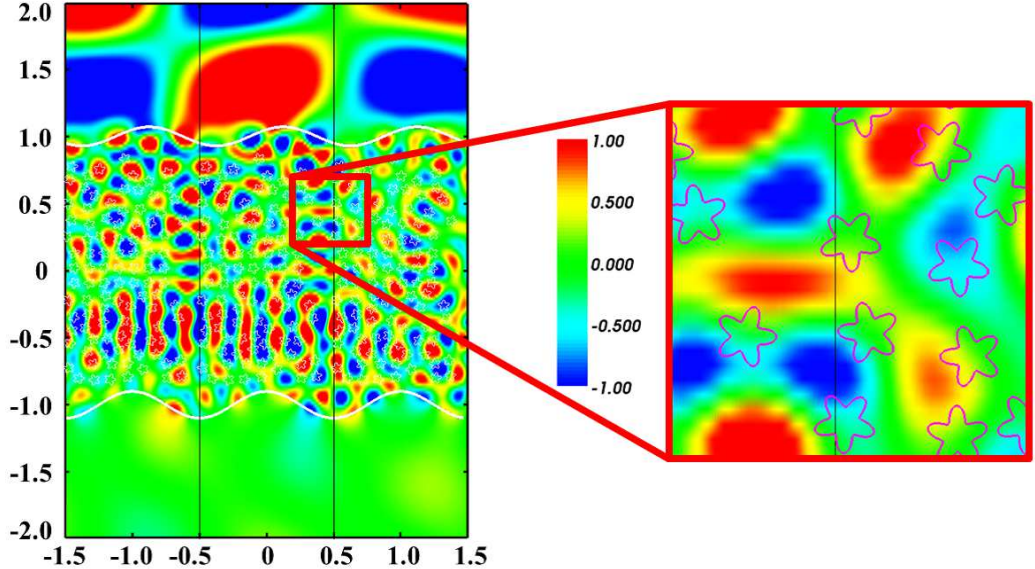}
\caption{Real part of the total field with $100$ dielectric inclusions, some of which touch the walls $L_2$ or $R_2$, randomly distributed in a three-layered medium (see Example 3).
The wavenumber for each particle is $k_p=8$, and the wavenumbers for the three layers are $k_1=5$, $k_2=30$ and $k_3=5$. 
Each inclusion is a smoothed five-pointed star, approximately $0.4$ wavelengths in size for $k_2$.}
\label{fig:NumericalResult3}
\end{center}
\end{figure}

\begin{table}[tbp] 
\begin{center}
\begin{tabular}{|c|c|c|c|c|}
\hline
$k_2$ & number of iterations & GMRES error & CPU time (sec.) & Flux error \\
\hline
$1$ & $13$ & $9.68$e-11 & $18.0$ & $8.54$e-10 \\
$10$ & $15$ & $1.39$e-11 & $18.9$ & $7.52$e-11 \\
$20$ & $43$ & $4.55$e-11 & $26.6$ & $4.38$e-10 \\
$30$ & $83$ & $5.12$e-11 & $37.8$ & $1.02$e-8 \\
\hline
\end{tabular}
\end{center}
\caption{Convergence behavior of GMRES, the CPU time required and flux error for $100$ inclusions embedded in the central layer, where $k_2$ is allowed to vary from $1$ to $30$ (see Example 3).}
\label{tab:NumericalResult3}
\end{table}

In our last example, we consider the scattering from $100$ smoothed pentagons ($a_3=5$) with $k_p=8$. In particular, we allow some of the inclusions intersect the wall $L_2$ and $R_2$, as shown in Fig. \ref{fig:NumericalResult3}. 
We set the thickness of the central layer to be $l=2.0$, and define parameters $a_1=0.0309,\ a_2=0.0103$ in Eq. \eqref{par1}. The boundary of the inclusion is discretized by $N=300$ equispaced points and the scattering matrix $S_p$ is truncated up to $p=10$. 
The interfaces $\Gamma_1$ and $\Gamma_2$ are given by the graphs
$y = 1 + 0.05\sin(2\pi x) + 0.05\cos(2\pi x)$
and $y = -1 + 0.1\cos(2\pi x)$ respectively.
The incident angle is set to be $\theta=-\arccos(1.0-2\pi /10.0)+0.1$, which is the same as example 1.

\begin{remark}
Careful readers might think that our method should fail in this case, since the multipole expansion for each inclusion is only valid outside the disk that encloses the inclusion, yet some target nodes on the walls lie
inside these disks.
However, it turns out our scheme is still valid, by design,
due to cancellation in the $C_{wp}$ block due to translational symmetry:
whatever field induced by $\Omega_p\in \Omega_2$ to $L_2$, it is equal to the field induced by $\Omega_p$ in the right copy of $\Omega_2$ to $R_2$ up to the Bloch phase. Once the subtraction is made in Eq. \eqref{bound2}, the two cancel each other.
Thus the only effects of particle multipole expansions on the discrepancy
between walls $L_2$ and $R_2$ are at a distance of around $Pd$ or more.
The walls are ``invisible'' to the particles in this scheme.
This is discussed more explicitly in \cite[Sec.~3.1]{BG2011} and \cite[Sec.~2.3]{mlqp}.
\label{cancel}
\end{remark}

The above observation has been confirmed by the numerical tests, as shown in \fref{fig:NumericalResult3} and Table~\ref{tab:NumericalResult3}.

Figure \ref{fig:NumericalResult3} shows the total field in the case $k_2=30$. We can see smooth field distortion due to the inclusions, even though some of the inclusions intersect the walls $L_2$ and $R_2$. In Table~\ref{tab:NumericalResult3}, we test higher-contrast materials: the wavenumber in the middle layer varies from $1$ to $30$ for the three-layered medium. In all cases, we can get $10$ digits of accuracy of GMRES and at least $8$ digits accuracy has been guaranteed in terms of flux error.

\section{Conclusion} 
\label{Conc}

We have demonstrated an efficient new scheme to solve the quasi-periodic boundary value problem arising when a time-harmonic plane wave is incident on a 
layered periodic structure containing a large number of inclusions,
as occurs in various composites and solar cell designs.
The method is based entirely on free-space Green's functions,
using an expanded linear system to
enforce quasi-periodicity and radiation conditions explicitly.
This avoids expensive computations of the quasi-periodic Green's function,
yet is robust at all scattering parameters including Wood's anomalies
(where the latter fails to exist).
We have shown high accuracies even at Wood's anomalies, and with inclusions
intersecting unit cell walls.

For the scattering between multiple inclusions, we introduce the scattering matrix for each inclusion and use it as a block-diagonal precondtioner, which greatly improves the conditioning of multi-particle scattering system. We also apply the FMM to accelerate the translation operator between the different structures. In the end, the system is solved iteratively by GMRES, scaling
optimally (linearly) in $M$ the number of inclusions, at fixed frequency.
This claim of $\bigO(M)$ scaling holds only if the number of iterations
is independent of $M$; however, our numerical tests suggest that any
growth with $M$ is very weak in our setting
(see Fig.~\ref{fig:NumericalResult1-1}).
Multi-particle systems with tens of thousands of unknowns are solved
to around 9-digit accuracy on a laptop in a few minutes.

There are several possible extensions that we leave for future work.
The method can easily be adapted to TE polarization, to complex permittivities,
and to non-smooth inclusions.
If higher aspect ratio unit cell regions (i.e.\ heights much bigger than the period) are needed, $P$ can be increased, although this
slows down the FMM which applies the large ${\mathcal D}$ matrix block.
High aspect ratios could instead be handled by replacing the
(intrinsically isotropic) $J$-expansions with proxy nodes
as in \cite{mlqp}, on oval curves.
Our scheme naturally generalizes to bi-periodic structures in 3D,
with the matrix $\mathcal A$
still directly invertible (with size of order $10^4$) at low frequencies.
Other future work includes a rigorous error analysis of the scheme.


\section{Acknowledgments}
The authors would like to thank Leslie Greengard for several useful discussions which greatly improved this paper.

\bibliographystyle{abbrv}
\bibliography{reference}

\end{document}